\documentclass{amsart}
\usepackage{amssymb,euscript,amsmath, mathrsfs}
\usepackage[dvips]{graphicx}
\usepackage[dvips]{color}

\newcounter{ENUM}
\newcommand{\itm}{\item}
\newenvironment{ilist}{\renewcommand{\theENUM}{\roman{ENUM}}\renewcommand{\itm}{\addtocounter{ENUM}{1}\item[(\theENUM)]}\begin{itemize}\setcounter{ENUM}{0}}{\end{itemize}}
\newenvironment{alist}[1][0]{\renewcommand{\theENUM}{\alph{ENUM}}\renewcommand{\itm}{\addtocounter{ENUM}{1}\item[\theENUM)]}\begin{itemize}\setcounter{ENUM}{#1}}{\end{itemize}}

\newcommand{\margh}[1]{}

\input xy
\xyoption{all}
\CompileMatrices

\def\P{{\mathbb P}}

\def\C{{\mathbb C}}

\def\ch{\operatorname{char}}

\def\gen{\operatorname{gen}}

\def\ram{\operatorname{ram}}
\def\branch{\operatorname{branch}}

\def\top{\operatorname{top}}

\def\tame{\operatorname{tame}}

\newtheorem{thm}{Theorem}[section]
\newtheorem{prop}[thm]{Proposition}
\newtheorem{lem}[thm]{Lemma}
\newtheorem{cor}[thm]{Corollary}

\theoremstyle{definition}
\newtheorem{defn}[thm]{Definition}

\newtheorem{ex}[thm]{Example}

\theoremstyle{remark}
\newtheorem{notn}[thm]{Notation}
\newtheorem{rem}[thm]{Remark}

\numberwithin{equation}{section}

\begin{document}
\title{Linear series and existence of branched covers}
\author{Brian Osserman}
\begin{abstract}
In this paper, we use the perspective of linear series, and in particular
results following from the degeneration tools of limit linear series, to
give a number of new results on existence and non-existence of tamely
branched
covers of the projective line in positive characteristic. Our results are
both in terms of ramification indices and the sharper invariant of monodromy
cycles, and the first class of results are obtained by intrinsically
algebraic and positive-characteristic arguments.
\end{abstract}
\thanks{This paper was supported by fellowships from the Japan Society for 
the Promotion of Science and the NSF}
\maketitle

\section{Introduction}

Over the complex numbers, the classical theory of branched covers and the 
Riemann existence theorem give a complete description of branched covers of 
curves in terms of, in the case of a base of genus $0$, the monodromy around 
branch points. Techniques of lifting
to characteristic 0 and comparing to the transcendental situation allow one 
to conclude that when (the order of) monodromy groups are prime 
to $p$, the situation in characteristic $p$ remains the same as the classical 
situation. However, when monodromy groups are divisible by $p$, even when
one studies tame covers the situation
becomes far more delicate. The main issue is that, although it remains true
that tame covers always lift to 
characteristic 0, it is no longer the case that a cover in characteristic 0
necessarily has good reduction to characteristic $p$. One subtlety in this
context is that existence and non-existence of tame covers can depend quite
strongly on the moduli of the curve; see Tamagawa's \cite{ta1}. We will 
avoid this 
issue entirely by focusing on the question of which tame covers exist for 
generic curves, but even in this context very little is known. 

Much of the work to date on existence of tame covers (see, e.g., Raynaud's 
\cite{ra3})
focuses on situations where one can still show that one has good reduction 
from characteristic 0. In this paper we pursue an entirely different tack, 
using degeneration techniques and the point of view of linear series. We 
therefore obtain results without lifting to characteristic 0 and invoking 
transcendental techniques, except that our results on monodromy groups are 
stated in terms of systems of generators which are only known to exist by 
transcendental methods. Although the idea of using degeneration techniques 
to obtain this sort of result is not new (see, e.g., Bouw-Wewers \cite{b-w} 
and Harbater-Stevenson \cite{h-s1}), the introduction of certain key new 
ingredients from the 
point of view of linear series allows for far stronger results than had 
previously been known, including in particular sharp non-existence results 
in the case of genus-0 covers of the projective line with either all 
ramification indices less than $p$, or only three ramification points. 

Due to the nature of the arguments, our results are all stated in terms of 
covers of the projective line with only one ramified point over each of 
$r$ general branch points. This means that the group-theoretic situation
is rather special, and in particular monodromy groups are always
cyclic, alternating, or symmetric groups. 
However, the restrictions on both genus and branching type can be relaxed 
to a substantial degree via a combination of standard techniques; see \S
\ref{s-discuss} below.

Finally, in \S \ref{s-example}, we examine two elementary examples which 
demonstrate two phenomena: 
first, unlike in the case of characteristic $0$, in positive characteristic 
the particular possibilities for local monodromy cycles of covers can 
depend in a strong sense on which generators of the fundamental group of 
the base are used to obtain them; second, degeneration techniques seem to
be ``unreasonably effective'' in a sense that will be made precise below. 

We now fix some terminology so that we can state our main theorem more 
precisely. 

\begin{defn} Let $X$ be an $r$-marked curve of genus $0$ over an 
algebraically closed field, with marked points $Q_1,\dots,Q_r$. 
We say that a tuple 
$(\gamma_1,\dots,\gamma_r)\in (\pi^{\tame}_1(X))^r$
is a {\bf local generating system} for $\pi^{\tame}_1(X)$ if:
\begin{ilist}
\itm the $\gamma_i$ generate $\pi^{\tame}_1(X)$; 
\itm $\gamma_1 \cdots \gamma_r =1$;
\itm each $\gamma_i$ is a generator of an inertia group at $Q_i$
(i.e., the algebraic equivalent of a small loop around $Q_i$). 
\end{ilist}
\end{defn}

When we write $\pi^{\tame}_1(X)$ for a marked curve $X$, we always 
mean the tame fundamental group of the curve with the marked points removed.

We recall \cite[Cor.\ XII.2.12]{sga1} that by specializing from topological 
generators in characteristic $0$, one sees that there always exists a local 
generating system for $\pi^{\tame}_1(X)$. 

\begin{defn} We say an $r$-tuple of cycles 
$(\sigma_1,\dots,\sigma_r) \in (S_d)^r$ is a {\bf Hurwitz factorization}
for $(d,r,\{e_1,\dots,e_r\})$ if:
\begin{ilist} 
\itm $\sigma_i$ has length $e_i$ for all $i$; 
\itm the product $\sigma_1 \cdots \sigma_r$ is trivial; 
\itm the $\sigma_i$ generate a transitive subgroup of $S_d$. 
\end{ilist}

We will also sometimes say simply that $(\sigma_1,\dots,\sigma_r)$ is 
a Hurwitz factorization when no $d$ or $e_i$ are specified, and conditions
(ii) and (iii) are satisfied.
\end{defn}

If we are given a tame cover $f$ of $X$, and a local generating system,
then we obtain a Hurwitz factorization $(\sigma_1,\dots,\sigma_r)$ by
labeling the fiber over the base point, and letting $\sigma_i$ be the 
monodromy of $f$ induced by $\gamma_i$.

Recall that the pure braid group acts on local generating systems: 
the $i$th generator of the braid group $B_r$ acts by replacing 
$(\gamma_i,\gamma_{i+1})$ by 
$(\gamma_{i+1},\gamma_{i+1}^{-1}\gamma_i \gamma_{i+1})$, but doesn't
respect the ordering. The pure braid group is the kernel of the natural
map $B_r \to S_r$, so respects the ordering and gives a well-defined action.
In the classical settings, any two topological local generating systems
are related by a pure braid transformation. 

In the case that $r=3$ or that $e_i<p$ for all $i$, will define in 
\S \ref{s-monodromy} below a purely group-theoretic condition
for a Hurwitz factorization to be {\bf $p$-admissible}.
We will see in Corollary \ref{num-comp} below
that $p$-admissibility is in fact a condition only on the $e_i$,
and can be roughly summarized as the requirement that the $e_i$ not
be too close to certain multiples of powers of $p$ in the case that
$r=3$, or that there be enough 
cycles of moderate length in the case that $e_i<p$ for all $i$
(see also Lemma \ref{3pt-num-lem} and Example \ref{non-ex-ex} below).

For the sake of comparison, we recall the genus-0 case of the following
theorem from {\it SGA}:

\begin{thm}\label{sga} \cite[Cor.\ XII.2.12]{sga1} Let $X$ be a
curve of genus $0$ over an algebraically closed field $k$, with
marked points $Q_1,\dots,Q_r$. Then there
exists a local generating system $(\gamma_1,\dots,\gamma_r)$ for 
$\pi_1^{\tame}(X)$ such that if $(\sigma_1,\dots,\sigma_r)$ is any Hurwitz 
factorization for $(d,r,\{e_1,\dots,e_r\})$ satisfying:
\begin{ilist}
\itm $2d-2=\sum_i(e_i-1)$; 
\itm either $\ch k=0$, or $\ch k =p>0$, and the group 
$\left<\sigma_1,\dots,\sigma_r\right> \subseteq S_d$ has order prime to $p$, 
\end{ilist}
there exists a cover $f:\P^1 \rightarrow X$ of degree $d$, with $f$ branched 
only over the $Q_i$, and such that for all $i$, the local monodromy around 
$\gamma_i$ is given by $\sigma_i$. 
\end{thm}

Our main theorem is then the following:

\begin{thm}\label{main} Let $(\sigma_1,\dots,\sigma_r)$ be a Hurwitz
factorization for $(d,r,\{e_1,\dots,e_r\})$, with $2d-2=\sum_i(e_i-1)$, 
every $e_i$ prime to $p$,
and where we suppose in addition that either $r=3$, or $e_i<p$ for all $i$.
Fix also a local generating 
system $(\gamma_1,\dots,\gamma_r)$ for $\pi_1^{\tame}(X^{\gen})$, where 
$X^{\gen}$ is the geometric generic $r$-marked curve of genus $0$, with 
marked points $Q_1,\dots,Q_r$.

Then the following are equivalent:
\begin{alist} 
\itm the tuple $(\sigma_1,\dots,\sigma_r)$ is $p$-admissible;
\itm there exists a map $f:\P^1 \rightarrow \P^1$ of degree $d$, and
distinct $P_1,\dots,P_r$ on the source $\P^1$, such that $f$ is ramified
to order $e_i$ at each $P_i$;
\itm there exists a cover $f:\P^1 \rightarrow X^{\gen}$ of degree $d$, and
a choice of local generating system $(\gamma_1',\dots,\gamma_r')$, with $f$
branched only over the $Q_i$, and such that for all $i$, 
the local monodromy around $\gamma'_i$ is given by $\sigma_i$; 
\itm there exists a cover $f:\P^1 \rightarrow X^{\gen}$ of degree $d$, and
a pure-braid transformation $(\gamma_1',\dots,\gamma_r')$ of 
$(\gamma_1,\dots,\gamma_r)$, with $f$ branched only over the $Q_i$, and 
such that for all $i$, the local monodromy around $\gamma'_i$ is given by 
$\sigma_i$; 
\end{alist}

Furthermore, if $r=3$, no pure braid transformation is required in d).
\end{thm}

See also Theorem \ref{three-pts} and Theorem \ref{numerical} below for a 
purely numerical criterion on
the $e_i$ determining whether b) holds. The necessity of working with
general branch points in order to give any purely group-theoretic criterion
for existence of covers with given monodromy is, as mentioned before, well
known, but we will give elementary examples in \S \ref{s-example} below.
The same examples will also justify the need for a pure-braid operation 
in d) of the Theorem as soon as $r>3$. Specifically, we will see:

\begin{prop}\label{non-invar} Let $\Sigma$ be a set of local generating 
systems for $\pi^{\tame}_1(X^{\gen})$ which is closed under pure-braid 
operations (e.g.,
the set of local generating systems arising as specializations of topological
systems). Then for $r=4$, and $p=3$ the total set of Hurwitz factorizations 
$(\sigma_1,\dots,\sigma_r)$ which arise as monodromy of tame covers (with
$d$ and $e_i$ allowed to vary arbitrarily) around
$(\gamma_1,\dots,\gamma_r) \in \Sigma$ depends on the choice
of $(\gamma_1,\dots,\gamma_r)$. 

In particular, there is no group-theoretic
criterion for recognizing when a Hurwitz factorization occurs as monodromy
of a tame cover around a local generating system which works simultaneously 
for all systems in $\Sigma$.
\end{prop}

The general thrust of the argument for our main theorem is to use the point
of view of linear series, and the associated degeneration tools of limit 
linear series,
to draw conclusions in the context of branched covers. There are four key 
ingredients: the 
exploitation of existence or non-existence of inseparable maps to draw 
conclusions about separable maps, as in \cite[Thm.\ 4.2]{os7}; the precise 
description \cite[Thm.\ 6.1]{os7} of when separable maps can degenerate to 
inseparable maps; a finiteness result \cite[Thm.\ 5.3]{os6}, 
proved by relating the maps in question to 
certain logarithmic connections with vanishing $p$-curvature on the 
projective line, and applying results of Mochizuki; and finally, a result
of Liu and the author \cite{o-l2} in the classical setting showing that 
in the situation we
study, all Hurwitz factorizations always lie in a single braid orbit. 
The first is used for the case of three ramification points, while the 
second and third are used to conclude sharp non-existence results in the 
case that all ramification indices are
less than $p$. The last result is used to go from numerical results in
terms of ramification indices to statements in terms of monodromy cycles.

\section*{Acknowledgements}

I would like to thank Akio Tamagawa for encouraging me to study the
question addressed in this paper, urging me to pursue the sharper results 
in terms of monodromy groups, and describing existing results; I am also 
grateful for his and Shinichi Mochizuki's patience in explaining over many 
conversations background theory on monodromy groups and admissible 
covers. Finally, I would like to thank David Harbater, Jakob Stix, Robert 
Guralnick, Michel Raynaud and Robin Hartshorne for their helpful 
conversations.

\section{Numerical results: the case of three points}\label{s-base}

In this section, we prove our first sharp result, in the case of
covers $\P^1 \to \P^1$ with only three ramification 
points. We use the techniques of \cite[Thm.\ 4.2]{os7} to generalize the 
result there, giving a complete answer in this case. Our results here,
as in \S \ref{s-num} below, are expressed numerically in terms of 
ramification indices; we translate into monodromy group statements in
\S \ref{s-monodromy} below.

As in 
{\it loc.\ cit.}, the main idea is to evaluate the existence of a separable 
map with the given ramification by examining the possibility of an 
inseparable linear series with the same ramification. 
We therefore begin with a brief review of linear series in the case of 
dimension $1$, which is far simpler than the general case. 
For the general definitions, see \cite{os8}. 

For our purposes, a {\bf linear series} of dimension $1$
and degree $d$ on a curve $C$ consists of a map $f:C \to \P^1$ of degree
$d' \leq d$, together with an effective divisor of {\bf base points} on $C$ of 
degree $d-d'$. However, we consider two maps to give the same linear series if
they are related by an automorphism of $\P^1$. Thus, whereas branched covers
are considered up to automorphism of the source, linear series are 
considered up to automorphism of the target.

We say that a linear series is {\bf separable} or {\bf inseparable} 
depending on whether the map $f$ is separable or inseparable. We say that 
a linear series is {\bf ramified} to order (at least) $e$ at a point $P$ 
if the sum of the ramification index of $f$ at $P$ and the number of
base points at $P$ is (at least) $e$.

A direct computation with the Hurwitz formula implies that if we have 
$P_1,\dots,P_r$ on $\P^1$, and
$e_1,\dots,e_r$ such that $2d-2=\sum_i(e_i-1)$, then we can only have
a separable linear series on $\P^1$ ramified to order $e_i$ at every $P_i$ 
if the associated map $f$ has degree $d$, so that we have no base points.
Thus, we will typically be working with linear series that have no base
points, which are simply maps to $\P^1$ up to automorphism of the image.
The exception is when we work with inseparable linear series, as in
the proof of Theorem \ref{three-pts} below.

We first introduce some useful notation and terminology:

\begin{notn} Given a positive integer $e$ prime to $p$, we write 
$\bar{e}^{[m,u]}:=\lceil \frac{e}{p^m} \rceil$, and 
$\bar{e}^{[m,d]}:=\lfloor \frac{e}{p^m} \rfloor$. Also write 
$e^{[m,u]}:=p^m \bar{e}^{[m,u]} -e$, and 
$e^{[m,d]}:=e- p^m \bar{e}^{[m,d]}$.
\end{notn}

\begin{defn}\label{num-p-admiss-3pt} Given positive integers $(e_1,e_2,e_3)$
prime to $p$, satisfying the triangle inequality, and with $e_1+e_2+e_3$ 
odd, we say that 
the triple $(e_1,e_2,e_3)$ is {\bf numerically $p$-admissible} if
for any $m>0$, and any $S \subseteq \{1,2,3\}$ such that: 
\begin{ilist}
\itm $p^m \leq d$; 
\itm $e_i>p^m$ for all $i \in S$; 
\itm $\sum_{i\in S}\bar{e}_i^{[m,d]}+\sum_{i \not\in S} \bar{e}_i^{[m,u]}$ 
is odd,
\end{ilist}
the following inequality is always satisfied:
\begin{equation}\label{3pt-eq1}
\sum_{i\in S}e_i^{[m,d]}+\sum_{i \not\in S} e_i^{[m,u]} \geq p^m.
\end{equation}
Here $d$ is the integer with $2d-2=\sum_i(e_i-1)$.
\end{defn}

Note that the triangle inequality condition on $(e_1,e_2,e_3)$, i.e.,
$e_1 \leq e_2+e_3, e_2 \leq e_1+e_3$, and $e_3 \leq e_1+e_2$, is 
equivalent to the condition that $e_i \leq d$ for all $i$.

The notion of numerical $p$-admissibility intuitively corresponds to not
having all three ramification indices too close to certain positive 
multiples of $p^m$, for any $m$. The following reformulation is also
useful.

\begin{lem}\label{3pt-num-lem} We may replace \eqref{3pt-eq1} in the above 
definition as follows:
for a given $m,S$, if we denote by $d^{[m,S]}$ the integer satisfying 
$$2d^{[m,S]}-2=
\sum_{i\in S}(\bar{e}_i^{[m,d]}-1)+\sum_{i \not\in S} (\bar{e}_i^{[m,u]}-1),$$
then \eqref{3pt-eq1} is equivalent to
\begin{equation}d<p^md^{[m,S]}+\sum_{i \in S} e_i^{[m,d]}.
\end{equation}
\end{lem}

\begin{proof}
The equivalence of the two conditions may be checked directly 
from the definitions, by verifying the identity 
$$d-p^m d^{[m,S]}-\sum_{i \in S} e_i^{[m,d]} 
= \frac{1}{2}(p^m-1-\sum_{i\in S}e_i^{[m,d]}-\sum_{i\not\in S} e_i^{[m,u]}).$$
\end{proof}

Our main result gives an explicit and purely numerical criterion for the 
existence of maps with three ramification points in positive characteristic.
The main idea is that an inseparable map will exist (and hence, a separable 
map won't exist) with the desired ramification only if the $e_i$ are 
too close to appropriate positive multiples of $p^m$, violating numerical
$p$-admissibility.

\begin{thm}\label{three-pts} Suppose we are given positive integers 
$d, e_1, e_2, e_3$ with $2d-2=\sum_i(e_i-1)$, each $e_i$ prime to $p$, and 
the $e_i \leq d$ for all $i$,
together with distinct points $Q_1, Q_2, Q_3$ on $\P^1$. Then the triple 
$(e_1,e_2,e_3)$ is numerically $p$-admissible if and only if
there exists a (necessarily unique up to automorphism) separable cover 
$f:\P^1 \rightarrow \P^1$ of degree $d$, branched over each $Q_i$ with
a single ramification point of index $e_i$.
\end{thm}

\begin{proof} We first observe that it doesn't matter whether we fix 
ramification points $P_i$ or branch points $Q_i$. Indeed, any three points 
on either the source or target are automorphism-equivalent, and when there 
are only three ramification points, we have $e_i+e_j>d$ for any $i,j$, so 
no two ramification points can lie above a single branch point.

Since all ramification is specified, existence of a separable map is
equivalent to existence of a separable linear series with the given
degree and ramification. By \cite[Thm.\ 4.2, (i)]{os7}, such a separable
linear series exists if and only if there does not exist an inseparable 
linear series with (at least) the specified ramification.

Now, suppose we have such an inseparable linear series, corresponding to a
map $f$ of degree $d''$ with $d-d''$ base points; we can then write 
$d''=p^m d'$, where $f$ is the composition of the $m$th power of 
(the relative) Frobenius with a separable map of degree $d'$. Write 
$e_i'$ for the ramification indices at the $P_i$ of this separable map. 
Let $S \subseteq \{1,2,3\}$ be the subset consisting of $i$ with 
$p^m e_i'<e_i$. Thus, for each $P_i$ in $S$, to have the required 
ramification we need to have at least $e_i-p^m e_i'$ base points in our 
inseparable linear series. Therefore, we must have 
$d \geq p^m d' + \sum_{i\in S} (e_i-p^m e_i')$, which we check directly may 
be rewritten as
$$\sum_{i\in S}(e_i-p^m e'_i)+\sum_{i\not\in S}(p^m e'_i - e_i)\leq p^m-1 
<p^m.$$ 
This implies that we must have $e_i- p^m e_i' = e_i^{[m,d]}$ for all 
$i \in S$, and $p^m e_i' - e_i = e_i^{[m,u]}$ for all $i \not\in S$. 
We then find that 
$$\sum_{i\in S}\bar{e}_i^{[m,d]}+\sum_{i \not\in S} \bar{e}_i^{[m,u]}=
\sum_i e'_i$$ 
is odd, 
and furthermore that Equation \ref{3pt-eq1} is violated, 
completing one direction of the proof.

Conversely, we may suppose that for some $m, S$, satisfying 
the conditions of Definition \ref{num-p-admiss-3pt}, we have 
$d\geq p^md^{[m,S]}+\sum_{i \in S} e_i^{[m,d]}$. It clearly suffices to 
show that there exists an $f'$ of degree $d^{[m,S]}$ having ramification 
at least $e'_i:=\bar{e}_i^{[m,d]}$ for $i\in S$ and at least 
$e'_i:=\bar{e}_i^{[m,u]}$ for $i\not\in S$, since we can then compose with 
the $m$th power of Frobenius and add $e_i^{[m,d]}$ base points at $P_i$ for 
each $i\in S$ to obtain an inseparable map of degree $d$ with ramification 
at least $e_i$. In fact, it is enough to have $f'$ any linear series 
(separable or inseparable, with or without base points), as long as it 
has the specified degree, and ramification $e'_i$ at each $P_i$. In this
case, we ``compose with the $m$th power of Frobenius'' by composing the
associated map to $\P^1$ with the $m$th power of Frobenius, and 
multiplying the base point divisor by $p^m$; this produces a new linear 
series, multiplying both the degree and all the ramification indices by 
$p^m$. 

Thus, as long as the $e'_i$ satisfy the triangle inequality (which we
recall is equivalent to the condition that $e'_i \leq d^{[m,S]}$ for all 
$i$), by \cite[Thm.\ 4.2, (i)]{os7}
we can find such an $f'$. Although the triangle 
inequality for the $e_i$ alone is not enough to show the triangle inequality 
for the $e'_i$, one can check using the equivalence of the two inequalities 
of the lemma that the triangle inequality for the $e_i$ together with the 
hypothesized inequality $d\geq p^md^{[m,S]}+\sum_{i \in S} e_i^{[m,d]}$ 
does in fact imply that the $e'_i$ satisfy the triangle inequality, giving 
us our inseparable map of degree $d$, and thereby showing that no 
separable map can exist.
\end{proof}

We also recall:

\begin{cor}\label{3pt-easy} In the situation of the theorem, if also 
$e_1,e_2$ are less than $p$, numerical $p$-admissibility is equivalent to 
the condition that $d<p$.
\end{cor}

\begin{proof} This may be derived directly from the theorem, but also 
predates it, see \cite[Thm.\ 4.2, (ii)]{os7}
\end{proof}

We provide some examples to demonstrate the usage 
of the combinatorial condition of the theorem.

\begin{ex}\label{3pt-ex} 
The indices $(1,d,d)$ for $d$ prime to $p$ provide 
trivial examples of ramification indices being far enough away from 
multiples of $p^m$ that an inseparable map never exists, as exhibited by 
the separable map $x^d$. Less trivial is the corollary, for instance, that 
if we take indices $(2, d-1, d)$ for $d \not\equiv 0,1 \pmod{p}$ and $p>2$, 
we necessarily obtain a separable map.

More substantively, we examine the case that $e_i<2p$ for all $i$. The above
corollary treats the case that at least two of the $e_i$ are less than $p$,
so we may assume that at most one $e_i$ is less than $p$. To rule out
exceptional cases, we also assume $p>3$. First, let us
suppose that $e_1<p$, but $e_2,e_3>p$. In this case, a map exists if and
only if $d<2p$, and $e_2+e_3-e_1 \geq 2p$. On the other other hand, if 
$e_i>p$ for all $i$, a map exists if and only if $d \geq 2p$, and
$e_2+e_3-e_1,e_1+e_3-e_2,e_1+e_2-e_3 \leq 2p$.

For both cases, $m=1$ is the only possibility. In the first case, we see
that the possibilities for $S$ are $S=\{2,3\}$ or $S=\emptyset$, which 
corresponds to the two inequalities. In the second case, we can have 
$S=\{1,2,3\},\{1\},\{2\}\text{ or }\{3\}$, corresponding to the four given
equalities. 

We see that for the first case, the requirement is that $d$ is
not too large, and the $e_i$ are not too symmetric, whereas in the second
case, the requirement is that $d$ is sufficiently large, and the $e_i$ are
sufficiently symmetric. This may seem counterintuitive, but in fact makes
sense: if $d \geq 2p$, we can compose the Frobenius map with a degree $2$ 
map to obtain a map ramified to order $2p$ at two points, and to order $p$ 
elsewhere; in particular, if $e_1<p$ and $e_2,e_3<2p$, then such a map 
immediately gives the required ramification. Similarly, the more symmetric
the $e_i$ are, the easier it is to obtain the ramification by starting with
the Frobenius map and adding base points. In contrast, once we have 
$e_1,e_2,e_3>p$, we need them to be large enough that we cannot simply add
base points to the Frobenius map, and symmetric enough that the map of degree 
$2p$ is too assymetric to yield the required ramification.
\end{ex}

\section{Branched covers and linear series}\label{s-prelim}

Before proceeding to treat the case of covers with more branch points,
we need to develop some preliminary results relating the existence 
questions for linear series and branched covers. That is, we compare the 
situation for existence of maps with given ramification points on a fixed 
source curve, to maps with branch points fixed on the base. The reason for 
this is two-fold: first, our degeneration techniques for linear series are 
best-suited to answer existence questions only for general configurations 
of ramification points, and we would like to be able to strengthen such 
results; and second, the relationship between the two perspectives has 
certain subtleties, as illustrated by the first remark below. 

Our main result is the following.

\begin{thm}\label{prelim} Let $d,r, e_1, \dots, e_r$ be positive integers 
with $e_i \leq d$ and prime to $p$ for all $i$, and 
$2d-2+2g = \sum _i (e_i-1)$. Then the following are equivalent:
\begin{alist}
\itm there exists a smooth curve $C$ of genus $g$ with distinct $P_i \in C$ 
and a separable map $f:C \rightarrow \P^1$ ramified to order $e_i$ at the 
$P_i$;
\itm for general $Q_i \in \P^1$, there exists a smooth curve $C$ of genus $g$
and a separable map $f:C \rightarrow \P^1$ of degree $d$, branched over 
each $Q_i$, with a single ramification point of order $e_i$.
\end{alist}

If further $g=0$ and all the $e_i$ are less than $p$, we also have the 
following equivalent condition:
\begin{alist}[2]
\itm for general $P_i \in \P^1$, there exists a separable map 
$f:\P^1 \rightarrow \P^1$ ramified at each $P_i$ to order $e_i$.
\end{alist}
\end{thm}

\begin{proof} 
The equivalence of a) and b) is well-known via deformation arguments,
but also follows immediately from \cite[Cor.\ 3.2]{os3}; the idea of the 
latter is to combine a classical codimension count from the point of view 
of linear series with easy deformation theory of branched 
covers to compute the total dimension of the space of maps from curves of
genus $g$ to $\P^1$ with the desired ramification. 

Now, we want to show that if $g=0$ and all $e_i$ are less than $p$, we also 
have that a) implies c). By \cite[Appendix]{os7}, we have a moduli scheme 
$MR:=MR(\P^1,\P^1, (e_1,\dots,e_r))$ parametrizing tuples 
$(f, (P_1, \dots, P_r))$, where $f:\P^1 \rightarrow \P^1$ is a separable 
morphism of degree $d$, the $P_i$ are distinct $k$-valued points of the 
domain, and $f$ is ramified to order at least $e_i$ at $P_i$. This comes 
with natural forgetful morphisms $\ram: MR \rightarrow (\P^1)^r$ and 
$\branch: MR \rightarrow (\P^1)^r$ giving the ramification points and 
branch points, which is to say the $P_i$ and $f(P_i)$ respectively. From 
the argument of \cite[Cor.\ 3.2]{os3}, we know that $MR$ has dimension 
exactly $r+3$ if it is non-empty.

The main tool for the argument is a finiteness result, which is currently 
only available in the situation that all the $e_i$ are odd, so we will 
induct on the number of even $e_i$, noting that it must always be an even 
number for the degree to be integral. In the base case that all the $e_i$ 
are odd, by \cite[Thm.\ 5.3]{os6} there can be only finitely many maps for 
any given choice of the $P_i$, up to automorphism of the image space. Thus, 
all the fibers of the $\ram$ morphism are at most $3$-dimensional, and 
since $MR$ has dimension $r+3$ if it is non-empty, we find that the $\ram$ 
morphism must dominate the $(\P^1)^r$ parametrizing ramification point 
configurations, as desired. For the induction step, without loss of 
generality we can suppose that $e_1, e_2$ are even. Suppose that $MR$ is 
non-empty and fails to dominate $(\P^1)^r$ under the $\ram$ morphism. 
$MR$ is necessarily dominant under the $\branch$ morphism, so let 
$(f,(P_1,\dots,P_r))$ be a point in it such the all the $f(P_i)$ are 
distinct. But then, by \cite[Lem.\ 5.2]{os7} we see that if we replace 
$e_1, e_2$ by $p-e_1, p-e_2$, and denote the new moduli scheme by $MR'$, 
then $\ram(MR')$ still fails to dominate $(\P^1)^r$, and at our specific 
$P_i$, we obtain a new $f'$ having ramification 
$p-e_1, p-e_2, e_3, \dots, e_r$ at the $P_i$, contradicting the induction 
hypothesis.
\end{proof}

\begin{rem}\label{pathology} Note that we cannot hope to drop the hypothesis 
that all $e_i$ are less than $p$ in order to obtain the equivalence of 
condition c). Indeed, for $p>2$ if we consider the family of functions 
$x^{p+2}+tx^p-x$
as $t$ is allowed to vary arbitrarily, we see that we obtain an infinite
(tamely ramified) family of maps for which the ramification points remain 
fixed, while the branch points (necessarily) move. While such maps do 
exist, they occur only for special configurations of ramification points; see 
\cite[Prop.\ 5.4]{os7} for details, and Example \ref{ex-tot-ram} below for a
more detailed examination of a similar example. Note, however, that 
(still in the case 
$g=0$) our argument shows that the equivalence of c) will hold whenever we 
know that there are only finitely many linear series with the given 
ramification for arbitrary distinct configurations of the $P_i$.
\end{rem}

\begin{rem} The finiteness result of \cite{os6} used to show the equivalence 
of condition c) is atypical in that it is proved via a relationship to 
certain connections with vanishing $p$-curvature, Mochizuki's dormant 
torally indigenous bundles. The fundamental obstruction to obtaining a 
direct proof is that in the setting of connections, one first enlarges to 
the category of connections with nilpotent $p$-curvature before proving 
finiteness, and it is not clear what the analogous construction would be in 
the context of rational functions on $\P^1$.
\end{rem}

\section{Numerical results: the case that $e_i<p$ for all $i$}\label{s-num}

We now put together the results of \S \ref{s-prelim} with the results
on linear series of \cite{os7} to obtain a sharp statement on the 
existence and non-existence of covers $\P^1 \to \P^1$ with given 
ramification indices, in the case that all indices are less than $p$. Our 
methods are algebraic and intrinsically positive characteristic.

We first state the numerical condition which will be equivalent to existence
of tame covers in the case that $e_i<p$.

\begin{defn}\label{num-p-admiss} Fix $e_1, \dots, e_r$ 
positive integers, with $\sum_i(e_i-1)$ even, and further assume that
$e_i < p$ for all $i$.

We say that 
$(e_1,\dots,e_r)$ is {\bf numerically $p$-admissible} if
there exist $r-3$ positive integers $e'_2, \dots, e'_{r-2}$, prime to 
$p$, such that for any $m$ with $1 \leq m <r-1$: 
\begin{ilist}
\itm the triple $(e'_m,e_{m+1},e'_{m+1})$ satisfies the triangle inequality;
\itm $e'_m+e_{m+1}+e'_{m+1}$ is odd and less than $2p$;
\end{ilist}
where we use the convention $e'_1:=e_1$ and $e'_{r-1}:=e_{r}$.
\end{defn}

Note that for $r=3$, it follows from Corollary \ref{3pt-easy} that the 
definition of numerical $p$-admissibility agrees with the definition 
already given in \S \ref{s-base}.

\begin{figure}
\centering
\input{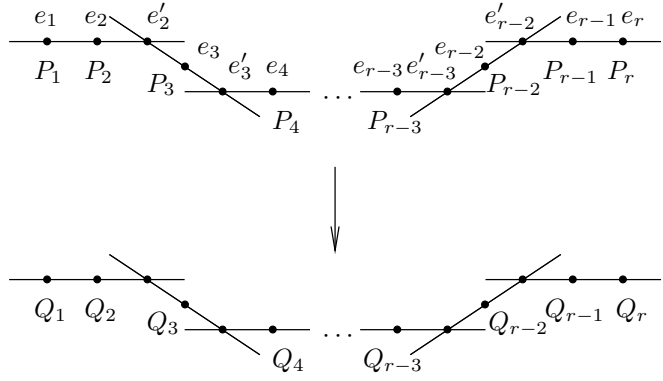}
\caption{The geometry behind numerical $p$-admissibility}
\label{fig-degen}
\end{figure}

In fact, 
numerical $p$-admissibility is equivalent to the existence of a
cover coming from smoothing a totally degenerate cover,
i.e., one constructed from a collection of maps $\P^1 \to \P^1$, 
with three ramification points each. The $e_i$ in this situation are the
ramification points above marked points, while the $e'_i$ are the 
ramification points above nodes; see Figure \ref{fig-degen}.

We further remark that there is a rather concrete interpretation of 
numerical $p$-admissibility when $e_i<p$ for all $i$, which 
says roughly that there have to be enough $e_i$ which aren't too large or 
small. See Example \ref{non-ex-ex} below for details.

Our theorem, which is a direct application of \S \ref{s-prelim}
and the results of \cite{os7}, is the following.

\begin{thm}\label{numerical} Fix $d, r$ and $e_1, \dots, e_r$ positive 
integers, with $e_i<p$ for all $i$ and satisfying $2d-2=\sum_i(e_i-1)$. 
Then the following are equivalent:

\begin{alist}
\itm The tuple $(e_1,\dots,e_r)$ is numerically $p$-admissible.
\itm For general $P_i \in \P^1$ there exists a separable map 
$f:\P^1 \rightarrow \P^1$ of degree $d$, ramified to order $e_i$ at each 
$P_i$.
\itm For general $Q_i \in \P^1$ there exists a separable map 
$f:\P^1 \rightarrow \P^1$ of degree $d$, branched over each $Q_i$ with a 
single ramification point of index $e_i$.
\itm There exists a separable map $f:\P^1 \rightarrow \P^1$ of degree $d$, 
ramified to order $e_i$ at some distinct $P_i$.
\end{alist}
\end{thm}

\begin{proof} The $g=0$ case of Theorem \ref{prelim} tells us that c)
and d) are equivalent, and that because all the $e_i$ are less than $p$,
we also have that b) and c) are equivalent. The equivalence of a) and b) 
is due to prior work using limit linear series and controlling degeneration 
from separable to inseparable maps; see \cite[Thm.\ 1.4]{os7}. 
\end{proof}

In order to illustrate how the (in principle rather complicated)
combinatorial conditions of the theorem may be applied, we examine some 
examples of existence and non-existence results.

\begin{ex}\label{num-ex} 
Existence examples are easy to construct: one simply needs to construct
a chain 
$$e_1,e_2,e'_2,\dots,e_{r-2},e'_{r-2},e_{r-1},e_r$$
such that each triple 
$(e_1,e_2,e_2'),(e_2',e_3,e_3'),\dots,(e'_{r-2},e_{r-1},e_r)$
satisfies numerical $p$-admissibility for three points. For instance,
since the triple $(p-2,\frac{p-1}{2},\frac{p-1}{2})$ is always
numerically $p$-admissible for $p>2$, if $r$ is a multiple of $3$ we can 
use the sequence
$$\frac{p-1}{2},p-2,\frac{p-1}{2},\frac{p-1}{2},p-2,\frac{p-1}{2},
\dots,\frac{p-1}{2},p-2,\frac{p-1}{2}$$
to obtain a map with $e_i=p-2$ when $i\equiv 2 \pmod{3}$, and 
$e_i=\frac{p-1}{2}$ otherwise (and we can make similar constructions in 
the cases that $r \equiv 1,2\pmod{3}$).
\end{ex}

For non-existence, the basic idea is that having $e_i$ very large or very 
small rigidifies the possibilities for the $e'_i$ (the former because 
of the condition $e'_m+e_{m+1}+e'_{m+1}<2p$, and the latter because of 
the triangle inequalities), and
can be used to ensure that no sequence of $e'_i$ can satisfy all conditions
at once. We will see that in essence, a map will exist unless all our 
$e_i$ are sufficiently close to $p$ or $1$, where ``sufficiently close''
depends on $p$ and on $r$. 

\begin{ex}\label{non-ex-ex}
For our first example, if $r$ is odd, and we set $e_i=p-1$ for $i<r-1$, and 
$e_{r-1}+e_r > p$ and odd, we see that the summation condition on 
$(e_1,e_2,e'_2)$ 
determines $e'_2=1$, and the triangle inequalities on $(e'_2,e_3,e'_3)$ 
determine $e'_3=p-1$, and so forth, alternating until $e'_{r-2}$ is 
determined as $p-1$. But then $(e'_{r-2},e_{r-1},e_r)$ violates the 
condition that their sum is less than $2p$.

Similarly, if $r$ is even, with $e_i=p-1$ for $i<r-1$, and 
$e_{r-1} \neq e_r$ with $e_{r-1}+e_r$ even, we find that $e'_{r-2}$ is
determined as $1$, and then
$(e'_{r-2},e_{r-1},e_r)$ violates the triangle inequality.

One can construct many more examples this way: if the $e_i$ are less than
$p-1$ but still close to $p$, there is more flexibility, so one has to
put restrictions on $r$. Similarly, one can use several small $e_i$, 
although as long as $e_i \geq 2$ each one will introduce some flexibility.
\end{ex}

Finally, we mention one other example of our results: the case of four points.

\begin{ex}\label{4pts-ex} Suppose that $r=4$, with $e_i<p$ for all $i$. 
We claim we have the
explicit formula that a cover exists if and only if $e_i > d+1-p$ for all 
$i$. Indeed, applying Theorem \ref{prelim} and \cite[Cor.\ 8.1]{os7} (which
is just an explicit computation with the combinatorial condition of the
corollary), we see that a map exists if and only if 
$$\min\{e_i,d+1-e_i,p-e_i,p-d-1+e_i\}_i>0.$$
We have $e_i,p-e_i,d+1-e_i>0$ by hypothesis, so the only possibility 
for non-existence is that
$p-d-1+e_i \leq 0$ for some $i$, which gives the desired statement.
\end{ex}

\begin{rem} Although the existence of maps with given ramification indices 
clearly does not depend on the order of the indices, our criterion is 
asymmetric. This reflects the fact that one can obtain the same result by 
degenerating to different totally degenerate curves, and obtain non-trivial 
combinatorial relations as a result. For details on some of these relations 
in a slightly different setting, see \cite{o-l1}.
\end{rem}

\section{Group-theoretic results}\label{s-monodromy}

In this section, we reformulate our previous numerical results on branched 
covers in terms of the sharper invariant of monodromy groups, ultimately
proving Theorem \ref{main}. Unlike the 
prior results, our results here will be dependent on transcendental 
techniques, as even our statements work systematically with local 
generating systems of tame fundamental groups. In \S \ref{s-discuss}
below, we discuss how the various restrictions of 
our main theorem can be relaxed in various ways to obtain large families 
of existence and non-existence results for additional covers.

The translation from numerical results to group-theoretic results is
much simpler in the three-point case, thanks to the following lemma.

\begin{lem}\label{3pt-hurw} Given $(d,3,\{e_1,e_2,e_3\})$, there exists a
unique Hurwitz factorization $(\sigma_1,\sigma_2,\sigma_3)$, up to 
simultaneous relabelling. 
\end{lem}

\begin{proof} One can argue the uniqueness using the discussion of
\S \ref{s-base}, but in fact the unique Hurwitz factorization in this case 
may also be described explicitly; see \cite[Lem.\ 2.1]{o-l2}.
\end{proof}

We next state the promised definition of $p$-admissibility. First,
note that the pure braid group acts on Hurwitz factorizations in exactly
the same manner as on local generating systems.

\begin{defn}\label{p-admiss} Let $(\sigma_1,\dots,\sigma_r)$ be a Hurwitz 
factorization for $(d,r,\{e_1,\dots,e_r\})$, where all $e_i$ are prime to
$p$, and $2d-2=\sum_i(e_i-1)$. 

If $r=3$, we say that $(\sigma_1,\sigma_2,\sigma_3)$ is {\bf $p$-admissible} 
if $(e_1,e_2,e_3)$ is numerically $p$-admissible (Definition 
\ref{num-p-admiss-3pt} above). 
 
If $e_i<p$ for all $i$, we say that $(\sigma_1,\dots,\sigma_r)$ is
{\bf $p$-admissible} if there exists a 
pure-braid transformation replacing
$(\sigma_1,\dots,\sigma_r)$ by $(\sigma'_1,\dots,\sigma'_r)$ and such that:
\begin{ilist}
\itm for any $m$ with $1 \leq m \leq r-1$, the partial product 
$\sigma''_m:=\sigma'_1 \cdots \sigma'_m$ is a cycle;
\itm for any $m$ with $1 \leq m < r-1$, the sum of the lengths of
$\sigma''_m,\sigma'_{m+1},\sigma''_{m+1}$ is less than $2p$.
\end{ilist}
\end{defn}

The condition for $r>3$ has a geometric interpretation in terms of totally 
degenerate covers; see 
the discussion following Definition \ref{num-p-admiss} as well as
Figure \ref{fig-degen} above.

As with numerical $p$-admissibility, if $r=3$ and $e_i<p$ for all $i$,
Corollary \ref{3pt-easy} implies that the two above definitions are 
equivalent (and indeed, that no pure braid transformation is necessary).

We also mention that if $(\sigma_1,\dots,\sigma_r)$ is $p$-admissible, 
then one checks easily that the tuple of lengths $(e_1,\dots,e_r)$ of the 
cycles is numerically $p$-admissible, by letting the $e'_i$ be the 
lengths of the $\sigma''_i$.

The basic lemma required to go from our numerical results to the final
group-theoretic statements is the following.
Although the proof of the lemma is by explicit construction, the intuition 
comes from considering monodromy groups of certain admissible covers. See 
also Figure \ref{fig-degen} above and Remark \ref{geom-rem} below.

\begin{lem}\label{exist-expl} Suppose that $e_i<p$ for all $i$, and 
$(e_1,\dots,e_r)$ is numerically 
$p$-admissible (Definition \ref{num-p-admiss}). Let $d$ be determined by
$2d-2=\sum_i(e_i-1)$. Then there exists a
Hurwitz factorization for $(d,r,\{e_1,\dots,e_r\})$ which is $p$-admissible.

More precisely, given any $e'_2,\dots,e'_{r-2}$ verifying numerical
$p$-admissibility, there is a Hurwitz factorization 
$(\sigma_1,\dots,\sigma_r)$ such that the partial products $\sigma''_i$
are cycles of length $e'_i$.
\end{lem}

\begin{proof} In fact, we produce a Hurwitz factorization which satisfies
the conditions for $p$-admissibility without any braid transformation. 
The proof is inductive, with a base case of $r=3$.  Indeed, the $r=3$ case 
is immediate from the definition, once we know that a Hurwitz factorization 
exists, which follows from Lemma \ref{3pt-hurw}.

For the induction step, suppose our assertion holds for $r-1$. Suppose also
we are given $(e_1,\dots,e_r)$, together with $e'_2,\dots,e'_{r-2}$ 
satisfying the conditions for the $r$-tuple to be numerically 
$p$-admissible. We then note that $(e_1,\dots,e_{r-2},e'_{r-2})$ is also
numerically $p$-admissible, with degree 
$d'=d-\frac{e_{r-1}+e_r-e'_{r-2}-1}{2}$. By the induction hypothesis, we
can find some Hurwitz factorization 
$(\sigma_1,\dots,\sigma_{r-2},\sigma'_{r-2})$ for 
$(d',r-1,\{e_1,\dots,e_{r-2},e'_{r-2})$ which satisfies the conditions for
$p$-admissibility without any braid transformation. We then note
that $(e'_{r-2},e_{r-1},e_r)$ is also numerically $p$-admissible,
with degree $d''=\frac{e'_{r-2}+e_{r-1}+e_r-1}{2}$. By the base case, we
can find a corresponding Hurwitz factorization 
$(\sigma''_{r-2},\sigma_{r-1},\sigma_r)$. In particular, the existence of
the two Hurwitz factorizations implies that we have
$e'_{r-2} \leq d',d''$. Now, we observe that
$d'+d''=d+e'_{r-2}$, so we see that $d',d'' \leq d$, and we can map 
$\{1,\dots,d'\}$ and 
$\{1,\dots,d''\}$ into $\{1,\dots,d\}$ such that:
\begin{ilist}
\itm $\{1,\dots,d'\}$ maps into $\{1,\dots,d'\}$ and $\{1,\dots,d''\}$
maps into $\{d'-e'_{r-2}+1,\dots,d\}$;
\itm $\sigma'_{r-2}$ maps to the inverse of $\sigma''_{r-2}$.
\end{ilist}
If we consider $(\sigma_1,\dots,\sigma_r)$ as lying in $S_d$ via these maps,
we then check easily that they give a Hurwitz factorization for 
$(d,r,\{e_1,\dots,e_r\})$, and furthermore satisfy the conditions for
$p$-admissibility without any braid transformation.
\end{proof}

We also recall the main theorem of \cite{o-l2}:

\begin{thm}\label{l-o-main} (Liu-Osserman) Given $(d,r,\{e_1,\dots,e_r\})$ 
with $2d-2=\sum_i(e_i-1)$, any two Hurwitz factorizations are related by 
a pure braid transformation.
\end{thm}

Putting together the theorem and the lemma (together with the earlier 
observation that $p$-admissibility implies numerical
$p$-admissibility), we conclude:

\begin{cor}\label{num-comp} A Hurwitz factorization 
$(\sigma_1,\dots,\sigma_r)$ for $(d,r,\{e_1,\dots,e_r\})$ is $p$-admissible 
if and only if $(e_1,\dots,e_r)$ is numerically $p$-admissible.
\end{cor}

We are now ready to prove Theorem \ref{main}.

\begin{proof}[Proof of Theorem \ref{main}] If 
$(\sigma_1,\dots,\sigma_r)$ is $p$-admissible, then $(e_1,\dots,e_r)$ is 
numerically $p$-admissible, so we see by Theorem \ref{three-pts} and 
Theorem \ref{numerical}
that there exists a cover $f$ with ramification indices 
$(e_1,\dots,e_r)$, and we conclude that a) implies b).

Now, a cover as in b) has monodromy $(\sigma'_1,\dots,\sigma'_r)$ 
around $(\gamma_1,\dots,\gamma_r)$ for some Hurwitz factorization
$(\sigma'_1,\dots,\sigma'_r)$. We then have by Theorem \ref{l-o-main}
that $(\sigma_1,\dots,\sigma_r)$ is related to $(\sigma'_1,\dots,\sigma'_r)$
by some pure braid transformation, so it follows that if we replace 
$(\gamma_1,\dots,\gamma_r)$ by $(\gamma'_1,\dots,\gamma'_r)$ under the same
transformation, the local monodromy of $f$ around 
$\gamma'_i$ is given by $\sigma_i$, so we have that b) implies d).

On the other hand, it is clear that d) implies c) which implies b), so
it only remains to check that b) implies a) under our hypotheses.
Accordingly, suppose we have a cover $f$ with ramification indices 
$(e_1,\dots,e_r)$ at points $P_1,\dots,P_r$, and suppose that either 
$e_i<p$ for all $i$, or that $r=3$. We claim that
$(e_1,\dots,e_r)$ is numerically $p$-admissible: indeed, in
the first case, this follows by Theorem \ref{numerical}, while in the
second, it follows from Theorem \ref{three-pts}.
Thus, by Corollary \ref{num-comp}, we see that 
$(\sigma_1,\dots,\sigma_r)$ is $p$-admissible, so in either
of these cases, we have b) implies a), as desired. 

The last assertion is
that for $r=3$, no braid operations are necessary: this follows from the
uniqueness of the Hurwitz factorization (Lemma \ref{3pt-hurw}) in that case.
\end{proof}

We give a simple example of the theorem.

\begin{ex} From the existence portion of our main 
theorem, we recover examples along the lines of covers constructed by
Bouw and Wewers \cite{b-w}, and Harbater and Stevenson \cite{h-s1}:

Suppose $p>3$ and $r \geq 3$, and fix general points $Q_1,\dots,Q_r$ on 
$\P^1$. Then there exists a tame cover $f:\P^1 \to \P^1$ of degree $d=r+1$, 
with monodromy group $A_d$, and with a single ramified point over each 
branch point, of order $3$.

By \cite[Thm.\ 5.3]{o-l2}, any such cover has monodromy group $A_d$, so
it suffices to note that the tuple 
$(\underbrace{3,\dots,3}_{r\text{ times}})$ is always numerically 
$p$-admissible for $p>3$ and $r \geq 3$. Indeed, if we set $e'_i=3$ for all 
$i$, we see that the conditions will be satisfied. 
\end{ex}

For some non-existence examples following from the theorem, see
\S \ref{s-discuss} below.

\begin{rem}\label{rem-imprim} 
Guralnick conjectures (see \cite[p.\ 2]{g-s1} for details) that other than
finitely many possibilities, all primitive genus-$0$ groups should be
cyclic, alternating, symmetric, or one of a few other possible families.
Thus, even though all the monodromy groups in Theorem \ref{main} are
one of the three types mentioned \cite[Thm.\ 5.3]{o-l2}, since they
are all primitive (except in the cyclic case with composite degree), the 
situation is actually reasonably general. That said, our non-existence 
results can also be immediately applied to draw conclusions on a far 
wider class of monodromy groups arising as imprimitive subgroups of $S_d$; 
see Theorem \ref{gen-non-ex} below.
\end{rem}

\section{Two examples}\label{s-example}

In this section we compute two elementary examples which shed light on the
subtleties of the situation in characteristic $p$, particularly regarding:
difficulty of using Riemann existence-style theorems to compute Hurwitz 
numbers; the peculiarly good behavior of degenerations; and (for the sake
of completeness, although this is well-known) the difficulties of dropping 
generality of branch points in existence statements. 
In particular, we prove Proposition \ref{non-invar}.

\begin{ex}\label{ex-simple-branch} We explicitly compute the situation
for covers of $\P^1$ of degree 3, with 4 simple branch points. We may assume
that $0$, $1$, and $\infty$ are three of the four ramification points, and
each is mapped to itself. We are thus considering functions of the form
$f(x)=\frac{ax^3+bx^2}{cx+d}$ where $a,b,c,d$ are all non-zero, and we can 
therefore choose to set $a=1$. Further, since $f(1)=1$, we have $d=b+1-c$.
We will also suppose that $f$ is simply ramified at some $\lambda$, and 
that $f(\lambda)=\mu$; we then want to compute $f$ in terms of $\mu$ to find 
all possibilities with particular prescribed branch points. Differentiating
$f$, and imposing the desired conditions, we find that
$c=\frac{2b^2+5b+3}{b+1}=2b+3$, and that $\lambda=\frac{-b^2-2b}{2b+3}$. 
Finally, the condition for $f(\lambda)=\mu$ can be simplified to
$$b^4+(2+8\mu)b^3+36\mu b^2+54 \mu b+27 \mu=0.$$

Over $\C$, this quartic corresponds to the four solutions one
obtains by classical cycle decompositions. Specifically, if we fix a local
generating system
$(\gamma_1,\dots,\gamma_4)$ for $\pi_1^{\top}(\P^1\smallsetminus\{P_i\}_i)$,
by considering monodromy around the $\gamma_i$, we have that covers are in 
one-to-one correspondence with Hurwitz factorizations for 
$(d=3,r=4,\{2,2,2,2\})$, which up to equivalence are:
\begin{gather*}
(1 2)(1 2)(2 3)(2 3)\\
(1 2)(2 3)(2 3)(1 2)\\
(1 2)(2 3)(3 1)(2 3)\\
(1 2)(2 3)(1 2)(3 1)
\end{gather*}

We now turn to the case of characteristic $3$. Here, the formula for $b$
in terms of $\mu$ reduces to
$$b^4+(-\mu-1)b^3=0,$$
and since $b$ must be non-zero, we obtain the unique solution $b=1+\mu$,
which corresponds to the function $f=\frac{x^2(x+1+\mu)}{(-\mu-1)x-\mu}$,
and which gives $\lambda=\mu$.
\end{ex}

\begin{ex}\label{ex-tot-ram} We next consider a similar example, with 
$d=4$, $e_1=4$, and 
$e_2=e_3=e_4=2$. Normalizing as before, we can write $f(x)=ax^4+bx^3+cx^2$
with $a,c$ non-zero, and $a+b+c=1$. Differentiating to impose ramification,
we find $b=4-2c$, and $\lambda=\frac{2c}{4(c-3)}$. Finally, imposing
$f(\lambda)=\mu$ gives us that $c$ is a root of
$$c^4-(4+16\mu)c^3+144\mu c^2-432\mu c+432 \mu=0.$$
The quartic polynomial corresponds to the Hurwitz factorizations
\begin{gather*}
(1 2 3 4)(1 2)(4 3)(3 1)\\
(1 2 3 4)(1 2)(1 4)(4 3)\\
(1 2 3 4)(1 2)(3 1)(1 4)\\
(1 2 3 4)(1 3)(1 4)(2 3)
\end{gather*}

As before, in characteristic $0$, (or with $p>3$) we have $4$ covers. But 
once again, in characteristic $3$, we see that we get a single cover, with
$c=\mu+1$. Note here the oddity that in characteristic $3$, we have 
$\lambda=-1$ is always fixed, so that although a cover exists for general
$\mu$, in fact a map does not exist for general $\lambda$, and there are
infinitely many for $\lambda=-1$. In particular, this map does not 
come from a totally degenerate map as in Figure \ref{fig-degen} of 
\S \ref{s-num} above. This is part of a much more general phenomenon;
see \cite[Prop.\ 5.4]{os7}.
\end{ex}

We wish to emphasize two phenomena. The first is the following: although it 
is a feature of the classical situation that, thanks to the Riemann 
existence theorem, which Hurwitz factorizations are realized as the 
monodromy of a cover around a local generating system is 
independent of the choice of system, the same statement fails in positive
characteristic, even for tame covers. Indeed, both these examples
demonstrate the phenomenon of Proposition \ref{non-invar} (and in particular
either example proves the proposition). In both cases, all four Hurwitz
factorizations lie in a single braid orbit, as can be checked directly or
follows from \cite{o-l2}. Since there is only a single cover in 
characteristic $3$, if we start with any local generating system, we see 
that which Hurwitz factorization we get varies between all four 
possibilities as we change the generating system by pure braid operations.

This phenomenon also makes it quite difficult to compute Hurwitz numbers
using a Riemann-existence-style theorem, as we see that although we can
identify which Hurwitz factorizations arise as the monodromy of covers
around some local generating system, there is no group-theoretic criterion
which works independent of the choice of local generating system, so it
seems quite subtle to use such criteria to actually count covers.

\begin{figure}
\centering
\input{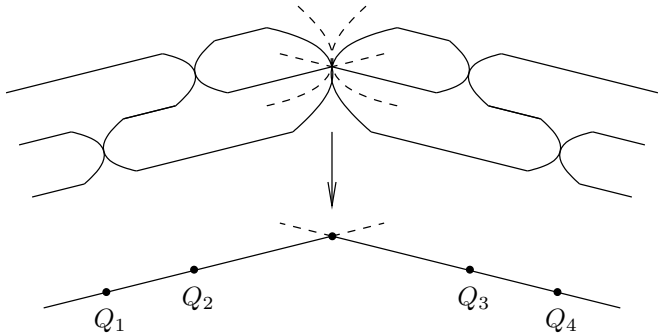}
\caption{The phantom admissible cover}\label{fig-admiss}
\end{figure}

The second phenomenon that we wish to discuss is that degenerations in
these examples behave ``unreasonably well'', in a different sense for each
example. In the case of the first example, we had originally hoped to give
an example where covers always degenerate from separable to inseparable.
Roughly, we start with a cover having monodromy cycles 
$(1,2),(2,3),(2,3),(1,2)$ at $P_1,\dots,P_4$ for some local generating
system (which we know exists). We then degenerate the base to a curve with 
two components, with $P_1,P_2$ on one component, and $P_3,P_4$ on the other. 
If the cover remained separable under this degeneration, it seems we would 
have to have a cover as in Figure \ref{fig-admiss}, 
with monodromy at the node given by $(3,2,1)$ and $(1,2,3)$ on the two
components, and this would not be possible in characteristic $3$. 
However, we see that we in fact cannot have bad degeneration: there is one 
cover in the smooth case, and also one cover (corresponding to monodromy 
$(1,2),(1,2),(2,3),(2,3)$, unramified over the node) after degeneration. 
What is going on is that for our analysis to work, the local generating 
system we choose on the generic fiber has to specialize to a ``geometric'' 
system in the sense of Remark \ref{geom-rem} below, so the only possible 
explanation is that while a cover with the specified monodromy exists for 
certain local generating systems, none of these systems specialize to 
geometric ones under the given degeneration.

On the other hand, in the second example, if we consider the same 
degeneration of the base, we see again that covers degenerate well, because 
again there is a single cover of both the smooth and degenerate curves.
In this case, the surprising aspect is that degeneration is well-behaved
from the point of view of covers, despite the fact that from the linear 
series point of view, the map cannot degenerate to a separable map (in 
fact, we saw that $\lambda$ cannot move at all). 
We thus see that in two different ways, situations where we might expect
to find a separable cover degenerating to an inseparable one do not in
fact give examples of bad degeneration.

We do however mention in this context that Bouw has an example, obtained
from Proposition 7.8 of \cite{bo2}, of a higher-genus cover with four 
branch point such that
degenerations are always inseparable. However, the argument in that case 
involves $p$-rank considerations, which one might hope to understand
completely. If there are no other obstructions to having good behavior
under (suitably generic) degeneration, it should be possible to use 
degeneration techniques to obtain a range of results beyond those 
presented here.

Finally, we mention that when $\mu=-1$ in either example,
the unique map specializes to an inseparable map.
We thus find that even though a separable cover exists for $\mu$ general,
we have a configuration of four distinct points, namely $-1,0,1,\infty$
over which no separable cover exists. Again, although this is already
well-known, we see that we have very simple examples showing the necessity
of considering general branch points in order to obtain any group-theoretic
condition for a Hurwitz factorization to be realized as monodromy of a 
cover.

\section{Further results and discussion}\label{s-discuss}

We begin with a discussion of the variety of ways in which the restrictions 
of Theorem \ref{main} can be relaxed. We state one generalized non-existence
result below, Theorem \ref{gen-non-ex}, and give examples. This is followed 
by a purely group-theoretic result which comes out of our arguments, and a
discussion of how one might use geometric arguments to prove a result like
Theorem \ref{main} without the benefit of the main result of \cite{o-l2}.

We first remark that standard deformation techniques can be used 
to vastly generalize our existence results. We don't state any general 
results here, because they tend to be complicated and we don't have any 
evidence that they are sharp. However, we do point out that many covers can
be constructed, where both the cover and the base may have higher genus,
by starting from the genus $0$ covers we have already constructed, 
gluing them together to obtain an admissible cover with higher genus, and
then deforming the admissible cover to obtain a smooth cover. 

Along these
lines, Fujiwara \cite{fu2} has used similar techniques to show that the 
prime-to-$p$ part of the fundamental group of any curve can be computed 
starting from the fundamental group of $\P^1$ with three marked points. 
One could broaden the array of examples even further by using the array of 
other covers known to existence, for instance from the papers \cite{ra3}, 
\cite{b-w} and \cite{h-s1} cited earlier, as well as additional work of
Stevenson \cite{st5}, Saidi \cite{sa1} and Raynaud \cite{ra4}.

Next, we observe that b) of Theorem \ref{main} makes no hypotheses
on the ramification points mapping to distinct branch points. Thus, 
by letting branch points coincide, and composing 
covers, we can considerably generalize our non-existence results:

\begin{thm}\label{gen-non-ex} Let $\{\sigma_1,\dots,\sigma_r\}$ be a
Hurwitz factorization of degree $d$ and genus $g$. That is, the $\sigma_i$
have trivial product and generate a transitive subgroup of $S_d$, and if 
$(e_1,\dots,e_m)$ are the lengths (with multiplicity) of all the cycles
in the disjoint cycle representations of the $\sigma_i$, we have
$2d-2+2g=\sum_i(e_i-1)$.

Suppose that the $\sigma_i$ act on blocks of size $m$,
so that they define permutations $\sigma'_1,\dots,\sigma'_r$ in $S_{d/m}$,
and suppose further that if $e_1,\dots,e_n$ are the lengths of all the
cycles in the disjoint cycle representations of the $\sigma'_i$, that:
\begin{ilist}
\itm either $e_i<p$ for all $i$, or $n=3$; 
\itm $2(d/m)-2=\sum_i(e_i-1)$;
\itm the tuple $(e_1,\dots,e_n)$ is not numerically $p$-admissible.
\end{ilist}

Then in characteristic $p$, there is no cover of degree $d$ and genus $g$ 
and with monodromy $\{\sigma_1,\dots,\sigma_r\}$ (around any choice of
local generating system).
\end{thm}

\begin{proof} Suppose that $\bar{f}$ is a cover with monodromy 
$\{\sigma_1,\dots,\sigma_r\}$ in characteristic $p$. Then we can lift
it to a cover $f$ in characteristic $0$ with the same monodromy. Under 
the stated hypotheses, we have that $f$ factors through a cover 
$g:\P^1 \to \P^1$ with monodromy $\{\sigma'_1,\dots,\sigma'_r\}$, and since
$f$ specializes to a separable cover in characteristic $p$, it follows 
that $g$ also specializes to a separable cover $\bar{g}$, having monodromy 
$\{\sigma'_1,\dots,\sigma'_r\}$. But the latter cover cannot exist by the
implication b) implies a) of Theorem \ref{main}. 
\end{proof}

We give two examples of the theorem, the first a genus-$0$ primitive cover,
and the second a higher-genus imprimitive cover.

\begin{ex} We see that the genus-$0$ cover of degree $9$, branched over 
three points with monodromy given by
$$\sigma_1=(1,2,3,4)(5,6,7,8), \sigma_2=(8,9,2,1)(4,3,6,5), 
\sigma_3=(1,5)(9,8,7,3),$$
doesn't exist in characteristic $p=5$, as the tuple $(4,4,4,4,4,2)$ is not
numerically $p$-admissible.
\end{ex}

\begin{ex} We consider the Hurwitz factorization
with $d=10$ and $g=1$ and three branch points, with monodromy
$$\sigma_1=(1,3,5,8,2,4,6,7), \sigma_2=(10,8,6,4,9,7,5,3), 
\sigma_3=(10,3,1,9,4,2)(7,8).$$
We see that this is imprimitive, acting on the blocks 
$[1,2],[3,4],[5,6],[7,8],[9,10]$ as cycles of length $4,4$, and $3$, so
this cover in characteristic $0$ factors through the genus $0$ cover of
degree $5$ corresponding to the Hurwitz factorization
$$\sigma_1'=(1,2,3,4),\sigma_2'=(5,4,3,2),\sigma_3'=(5,2,1).$$
In characteristic $5$, we have that the latter cover doesn't exist, so
we conclude that the genus-$1$ cover with the given monodromy also does
not exist in characteristic $5$.
\end{ex}

We next make the following simple group-theoretic observation
which we don't believe is obvious without some type of geometric argument: 

\begin{prop} Given a Hurwitz factorization $(\sigma_1,\dots,\sigma_r)$ for 
$(d,r,\{e_1,\dots,e_r\})$, where $2d-2=\sum_i(e_i-1)$, then there exists
a pure braid transformation $(\sigma'_1,\dots,\sigma'_r)$ of 
$(\sigma_1,\dots,\sigma_r)$ such that each partial product
$\prod_{i=1}^m \sigma'_i$ is a cycle, for $1 \leq m \leq r$.
\end{prop}

\begin{proof}
Indeed, if we start with such a Hurwitz factorization, if we consider
$p$ sufficiently large, the tuple is automatically $p$-admissible, so applying
Lemma \ref{exist-expl} we know there exists some $(\sigma'_1,\dots,\sigma'_r)$
of the desired form. Then by Theorem \ref{l-o-main} the two Hurwitz
factorizations are related by a pure braid operation. 
\end{proof}

In fact, the argument for Lemma \ref{exist-expl} shows that any $p>d$
is sufficiently large for the above argument. Of course, considering large 
$p$ above is equivalent to simply thinking about the 
situation in characteristic $0$. Geometrically, it might be simpler to
argue that we can switch to the linear series point of
view, degenerate to obtain a totally degenerate limit linear series, and
then switch back to the point of view of admissible covers. This ensures
that there is only one ramification point over each node, and hence that 
the partial products are cycles. This avoids any use of Theorem 
\ref{l-o-main}, but requires a good understanding of local generating
systems for admissible fundamental groups, discussed in more detail in
the remark that follows.

\begin{rem}\label{geom-rem} Because of Theorem \ref{l-o-main}, we have
been able to considerably simplify the transition from numerical to
group-theoretic results, compared to the arguments originally envisioned. 
However, the original arguments should hold more generally than 
Theorem \ref{l-o-main}, so we briefly sketch here how we would expect them
to go. We suppose we have a given branched cover, with given local 
generating system on the base. We wish to be able to say something about 
the monodromy of the cover, perhaps after pure braid transformation; for 
these purposes, it is probably best to assume (as is done in \cite{sga1}) 
that the local generating system arises as the specialization of a set of 
topological generators from characteristic $0$. The program is then as 
follows:

For the most general step, we suppose we have a family of semistable curves,
and a local generating system on the (smooth) geometric generic fiber,
arising from specialization from topological generators. We say that
a local generating system for the admissible fundamental group of the
special fiber is ``geometric'' if it arises by gluing local 
generating systems on each component of the normalization; in this case, 
one can express monodromy of admissible covers in terms of the monodromy
of components of the normalization. The basic assertion is that after
a pure braid transformation, the given local generating system on the 
geometric generic fiber specializes to a geometric local generating system
on the special fiber. Given appropriate definitions, the statement is clear 
in the topological setting, since any two choices of local generating 
systems on the smooth fiber are related by pure braid transformations. One 
then has to check that the definitions behave well with respect to 
algebraization and specialization. Finally, one would lift the original 
family of curves to characteristic $0$, and compare to the topological 
setting to obtain to desired result.

We return to considering a given branched cover of $\P^1$. We first think
of the cover as a $1$-dimensional linear series, and show that it can be 
degenerated to a limit linear series on some degeneration of the original 
curve. We then show that a family of limit linear series can be realized 
geometrically as a family of admissible covers. Applying the previous
discussion, after pure braid transformation the local generating system
specializes to a geometric local generating system, so we can use the 
geometry of the admissible cover to describe the monodromy of the smooth
cover.
\end{rem}

\bibliographystyle{hamsplain}
\bibliography{hgen}

\newcommand{\noopsort}[1]{} \newcommand{\printfirst}[2]{#1}
  \newcommand{\singleletter}[1]{#1} \newcommand{\switchargs}[2]{#2#1}
\providecommand{\bysame}{\leavevmode\hbox to3em{\hrulefill}\thinspace}
\begin{thebibliography}{10}

\bibitem{bo2}
Irene Bouw, \emph{The $p$-rank of ramified covers of curves}, Compositio
  Mathematica \textbf{126} (2001), no.~3, 295--322.

\bibitem{b-w}
Irene Bouw and Stefan Wewers, \emph{Alternating groups as monodromy groups in
  positive characteristic}, Pacific Journal of Mathematics \textbf{222} (2005),
  no.~1, 185--200.

\bibitem{fu2}
Kazuhiro Fujiwara, \emph{Etale topology and the philosophy of log}, Algebraic
  Geometry Symposium (Kinosaki), 1990, in Japanese, pp.~116--123.

\bibitem{sga1}
Alexandre Grothendieck, \emph{Revetements etales et groupe fondamental}, SGA,
  no.~1, Spring-Verlag, 1971.

\bibitem{g-s1}
Robert Guralnick and John Shareshian, \emph{Symmetric and alternating groups as
  monodromy groups of {R}iemann surfaces {I}: generic covers and covers with
  many branch points}, Memoirs of the AMS, to appear.

\bibitem{h-s1}
David Harbater and Katherine~F. Stevenson, \emph{Patching and thickening
  problems}, Journal Algebra \textbf{212} (1999), no.~1, 272--304.

\bibitem{o-l2}
Fu~Liu and Brian Osserman, \emph{The irreducibility of certain pure-cycle
  {H}urwitz spaces}, \mbox{arXiv:math.AG/0609118}.

\bibitem{o-l1}
\bysame, \emph{Mochizuki's indigenous bundles and {E}hrhart polynomials},
  Journal of Algebraic Combinatorics \textbf{23} (2006), no.~2, 125--136.

\bibitem{os8}
Brian Osserman, \emph{A limit linear series moduli scheme}, Annales de
  l'Institut Fourier, \mbox{arXiv:math.AG/0407496}, to appear.

\bibitem{os6}
\bysame, \emph{Mochizuki's crys-stable bundles: A lexicon and applications},
  Publications of {RIMS}, \mbox{arXiv:math.AG/0410323}, to appear.

\bibitem{os3}
\bysame, \emph{Deformations of covers, {B}rill-{N}oether theory, and wild
  ramification}, Mathematical Research Letters \textbf{12} (2005), no.~4,
  483--491.

\bibitem{os7}
\bysame, \emph{Rational functions with given ramification in characteristic
  $p$}, Compositio Mathematica \textbf{142} (2006), no.~2, 433--450,
  \mbox{arXiv:math.AG/0407445}.

\bibitem{ra3}
Michel Raynaud, \emph{Sp\'ecialisations des rev\^etements en caract\'eristique
  $p>0$}, Ann. Sci. \'Ecole Norm. Sup. \textbf{32} (1999), no.~1, 87--126.

\bibitem{ra4}
\bysame, \emph{Rev\^etements des courbes en caract\'eristique $p>0$ et
  ordinarit\'e}, Compositio Mathematica \textbf{123} (2000), no.~1, 73--88.

\bibitem{sa1}
Mohamed Saidi, \emph{Rev\^etements mod\'er\'es et groupes fondamental de
  graphe}, Compositio Mathematica \textbf{107} (1007), no.~3, 319--338.

\bibitem{st5}
Katherine~F. Stevenson, \emph{Galois groups of unramified covers of projective
  curves in characteristic $p$}, Journal of Algebra \textbf{182} (1996), no.~3,
  770--804.

\bibitem{ta1}
Akio Tamagawa, \emph{Finiteness of isomorphism classes of curves in positive
  characteristic with prescribed fundamental groups}, Journal of Algebraic
  Geometry \textbf{13} (2004), no.~4, 675--724.

\end{thebibliography}
\end{document}